\documentclass[12pt]{article}

\usepackage{amsmath,amssymb,amsfonts,theorem,makeidx,latexsym,epsfig,,subfigure}

\newtheorem{defn}{Definition}[section]

\newtheorem{lemma}[defn]{Lemma}

{\theorembodyfont{\rmfamily}

\newtheorem{ex}[defn]{Example}}

\newtheorem{thm}[defn]{Theorem}

\newtheorem{prop}[defn]{Proposition}

\newtheorem{cor}[defn]{Corollary}

\numberwithin{equation}{section}

\newcommand{\ltr}{ L^2(\mathbb R) }

\newcommand{\mn}{\mathbb N}

\newcommand{\mr}{\mathbb R}

\newcommand{\mz}{\mathbb Z}

\newcommand{\mc}{\mathbb C}

\newcommand{\mts}{ \{E_{mb}T_{na}g \}_{m,n \in \mz}}

\def\bp{{\noindent\bf Proof. \ }}

\def\ep{\hfill$\square$\par\bigskip}

\def\bqs{\begin{equation}}

\def\eqs{\tag*{$\square$}\end{equation}\par\bigskip}

\def\la{\langle}

\def\ra{\rangle}

\def\bop{\begin{op}\rm}

\def\eop{\end{op}}

\def\bee{\begin{eqnarray}}

\def\ene{\end{eqnarray}}

\def\bes{\begin{eqnarray*}}

\def\ens{\end{eqnarray*}}

\def\bei{\begin{itemize}}

\def\eni{\end{itemize}}

\def\bt{\begin{thm}}

\def\et{\end{thm}}

\def\bc{\begin{cor}}

\def\ec{\end{cor}}

\def\bpr{\begin{prop}}

\def\epr{\end{prop}}

\def\bl{\begin{lemma}}

\def\el{\end{lemma}}

\def\bd{\begin{defn}}

\def\ed{\end{defn}}

\def\bex{\begin{ex}}

\def\enx{\end{ex}}

\def\bfi{\begin{fig}}

\def\efi{\end{fig}}

\newcommand{\nft}{ || f||^2}

\title{On entire functions restricted to intervals, partition of unities, and dual Gabor frames}

\date{\today}

\author{Ole Christensen,
Hong Oh Kim,
Rae Young Kim}

\begin{document}

\maketitle

\begin{abstract} Partition of unities appear in many places in analysis.
Typically they are generated by compactly supported functions with a certain regularity. In this paper we consider partition of unities obtained as integer-translates of entire functions restricted to finite intervals. We characterize the entire functions that lead to a partition of unity in this way, and we provide characterizations of the  ``cut-off" entire functions, considered as functions of a real variable, to have desired regularity. In particular we obtain partition of unities generated by functions with small support and desired regularity.  Applied to Gabor analysis this leads to constructions of dual pairs of Gabor frames with low redundancy, generated by trigonometric polynomials with small support and desired regularity.
\end{abstract}

\begin{minipage}{120mm}

{\bf Keywords}\ {Entire functions, trigonometric polynomials, partition of unity, dual frame pairs, Gabor systems, tight frames}\\
{\bf 2000 Mathematics Subject Classification:} 42C40 \\

\end{minipage}
\

\section{Introduction} \label{187a}

Partition of unity conditions play an important role in many parts of mathematics, for example, within applied harmonic analysis \cite{Fe2, Fe4, FG2}. Usually the size of the support of the functions in the partition of unity and their regularity are key issues. Some of the most important partition of unities are obtained by considering integer-translates of an appropriately chosen function, e.g., a B-spline or a scaling function in the context of wavelet analysis.
In this paper we will consider  entire functions $P: \mc \to \mc$ which, for some fixed $N\in \mn,$ satisfy that
\bee \label{157af} \sum_{n\in \mz} P(x+n) \chi_{[0,N]}(x+n)=1,  \ x\in \mr.\ene The condition \eqref{157af} means exactly that the function $p:=P\chi_{[0,N]}$ satisfies the partition of unity condition. We will characterize the entire functions $P$ that are solutions to \eqref{157af}; in particular, we will see that such functions  automatically are $N$-periodic, i.e., they can be expanded in an everywhere convergent Fourier series. The Fourier series naturally connects with trigonometric polynomials. We will show that, still for any fixed $N\in \mn,$ given any $L\in \mn$ there exists a trigonometric polynomial $P$ such that the function $P\chi_{[0,N]},$ considered as a function of a real variable, belongs to $C^{2L-1}(\mr)$ and generates a partition of unity.

As an application of the results we will construct  dual  pairs of Gabor frames $\{E_{mb}T_{n}g\}_{m,n\in \mz}, \{E_{mb}T_{n}h\}_{m,n\in \mz}$ for
functions $g,h$ of  the form
$g  =  G \, \chi_{[0, N]}$ and $h =  H\, \chi_{[0, N]}$ for some trigonometric polynomials $G, H.$ In particular, our results show that for any $b\in ]0, 1/N]$ such constructions are possible, with $G$ and $H$ being trigonometric polynomials and the associated windows $g,h$ having  desired regularity; in contrast to the results in the literature, this is even possible for $N=2,$ i.e., higher smoothness is not obtained at the cost of larger support. Taking small values for $N,$ we obtain frames with low redundancy, generated by
functions with small support and desired regularity.
As a special case we obtain simple constructions where $G$ and $H$ are just powers of the sine function.
Finally, we  show that the condition $b\le 1/N$ is necessary  for such frame constructions to exist.

The paper is organized as follows. In the rest of this introduction we use standard frame theory to motivate the interest in the partition of unity condition \eqref{157af}. In Section \ref{157f} we carry out the analysis of this condition in a general fashion. Section \ref{277a} specializes to the case of trigonometric polynomials, which is a convenient setting for applications, and finally we connect with the Gabor analysis in Section \ref{157g}.

Gabor systems play a central role in time-frequency analysis. The basic idea is to decompose signals or functions into superpositions of certain time-frequency shifts of a fixed function $g.$ In the discrete case these time-frequency shifts have the form $\{e^{2\pi i mbx}g(x-na)\}_{m,n\in \mz}$ for appropriate parameters $a,b>0;$
using the translation operators $T_af(x):=f(x-a)$ and the modulation operators $E_bf(x):=e^{2\pi ibx}f(x),$ the time-frequency-shifts have the form of a coherent system $\mts.$ The system $\mts$ is called a {\it Gabor frame} if there exist constants $A,B>0$ such that
\bes A\, \nft \le \sum_{m,n\in \mz} | \la f, E_{mb}T_{na}g\ra|^2 \le B\, \nft, \ \forall f\in \ltr.\ens The system $\mts$ is a Bessel sequence if at least the upper frame condition is satisfied. It is well known that if $\mts$ is a Gabor frame, there exists at least one dual Gabor frame $\{E_{mb}T_{na}h\}_{m,n\in \mz},$ i.e., a Gabor frame such that the decomposition
\bes f= \sum_{m,n\in \mz} \la f, E_{mb}T_{na}g\ra E_{mb}T_{na}h\ens holds for all $f\in \ltr.$ The duality conditions by Ron \& Shen \cite{RoSh}, resp. Janssen
\cite{J} states that two Bessel sequences $\{E_{mb}T_{n}g\}_{m,n\in
\mz}$ and $\{E_{mb}T_{n}h\}_{m,n\in \mz}$ form dual frames for
$\ltr$ if and only if
\bee \label{gfs} \sum_{k\in \mz}
\overline{g(x+n/b+k)}h(x+k) & = & b\delta_{n,0}, \ a.e. \ x\in
\mr.\ene

For a bounded and compactly supported functions $g,$ the associated Gabor system $\mts$ automatically form a Bessel sequence. Furthermore, for functions $g,h$ having support in an interval of length $N,$ the condition \eqref{gfs} is automatically satisfied for $n\neq 0$ if we assume that $b\in ]0, 1/N].$
Thus, for functions $g  =  G \, \chi_{[0, N]}$ and $h =  H \, \chi_{[0, N}]$ as described in the introduction and for
$b\in ]0, 1/N],$ the Gabor systems
$\{E_{mb}T_{n}g\}_{m,n\in
\mz}$ and $\{E_{mb}T_{n}h\}_{m,n\in \mz}$ form dual frames if and only if the function
$P:=GH$ satisfies the condition
$\sum_{n\in \mz} P(x+n) \chi_{[0,N]}(x+n)=b.$
Discarding the factor $b$ then leads to the partition of unity constraint \eqref{157af}.

We will see that the connection to entire functions will bring some new aspects into the analysis.
Note that the connections between complex analysis and frame theory has already proved to be useful in other contexts, see, e.g., \cite{Ly,Se2,SW}. For more information on Gabor systems and frames we refer to the books \cite{G2, CBN}.

\section{Partition of unity for entire functions} \label{157f}

Motivated by the introduction we will consider entire functions
$P: \mc \to \mc$ satisfying the partition of unity condition \eqref{157af}.
We will first show that for such functions $P$ the restriction to $\mr$ in $N$-periodic. This implies that we have an extra tool at our disposal, namely, Fourier expansions.

\bl \label{157d} Let $N\in \mn.$ Then an entire function $P$  satisfies \eqref{157af} if and only if its restriction to $\mr$ is $N$-periodic and the Fourier coefficients $c_k$ in the expansion
\bee \label{157c} P(x)= \sum_{k\in \mz} c_ke^{2\pi ikx/N}, \ x\in \mr,\ene
satisfy that $c_k=  \frac1{N} \delta_{k,0}$ for $k\in N\mz.$\el

\bp Assume first that \eqref{157af} holds. Then, for $x\in [0,1],$
\bee \label{108a} P(x)+ P(x+1) + \cdots + P(x+N-1)=1.\ene
Since $P$  is an entire function, \eqref{108a} then holds for all $x\in \mr.$ Doing the similar calculation with  $x$ replaced by $x+1$ and subtracting the two expressions shows that
$P(x+N)=P(x), \ \forall x\in [0,1].$ The same calculation works with $[0,1]$ replaced by any interval $[n,n+1],$ so we conclude that the restriction of
$P$ to $\mr$ is $N$-periodic. Writing $P$ as the Fourier series \eqref{157c},
the equation \eqref{108a} takes the form
\begin{equation} \label{et-65}
 \sum_{k\in \mz} c_k\,\left[1+ e^{2\pi i k/N}
+ \cdots + \left( e^{2\pi i k/N}\right)^{N-1} \right] \, e^{2\pi ikx/N}=1.
\end{equation}
We note that
\begin{equation} \label{et-66}
1+ e^{2\pi i k/N}
+ \cdots + \left( e^{2\pi i k/N}\right)^{N-1}=\left\{
\begin{array}{ll}
   N, & k\in N\mz \\
   0, & k\notin N\mz.
\end{array}
\right.
\end{equation}
From \eqref{et-65} and \eqref{et-66}, we see that
$c_k=\frac{1}{N}\delta_{k,0}$ for $k\in N\mz.$
Conversely, if $P$ is $N$-periodic and  satisfies that
$c_k= \frac1{N}\, \delta_{k,0}$ for $k\in N\mz,$ then for $x\in [0,1]$,
\bes & \ & \sum_{n\in \mz} P(x+n) \chi_{[0,N]}(x+n) =
\sum_{n=0}^{N-1} P(x+n) \\ & = & \sum_{k\in \mz} c_k\,\left[1+ e^{2\pi i k/N}
+ \cdots + \left( e^{2\pi i k/N}\right)^{N-1} \right] \, e^{2\pi ikx/N}
= 1 \ens
by \eqref{et-66}. By periodicity \eqref{157af} holds for all $x\in \mr.$
\ep

It is well known that if an entire function is periodic when restricted to the real line, the Fourier coefficients $c_k$ have  exponential decay of arbitrary order, i.e., for any $a>0$ there exist a constant C such that $|c_k| \le C e^{-a|k|}$ for all $k\in \mz.$.Therefore, the expansion in \eqref{157c} converges even for all complex numbers $x$ and so the entire function P is $N$-periodic in the whole complex plane.

We will now fix $N\in
\mn$ and return to the solutions $P$ of \eqref{157af}. Our purpose is to investigate the
regularity of the functions $P\chi_{[0,N]}.$

\bt \label{et-68}
 Let $N\in \mn$. Assume that $P$ is an
$N$-periodic entire function
satisfying that $c_k=\frac{1}{N}\delta_{0,k}, k\in N\mz,$ and that the restriction of $P$ to $\mr$ is real-valued.  Then the following hold.
\begin{itemize}
\item[\rm{(a)}] There does not exist $P$ of this form
such that $P  \chi_{[0, N]}\in C^{\infty}(\mr)$; \item[\rm{(b)}] Fix $L \in \mn.$
Then  $P  \chi_{[0, N]} \in C^{L-1}(\mr)$ if and only if
\begin{equation} \label{et-78}
P(x)= \left(e^{\pi i x/N} \sin(\pi
x/N)\right)^{L} A_{L}(x)
\end{equation}
 for an  $N$-periodic entire function
$A_{L}(x):= \sum_{k \in \mz} a_k e^{2\pi i k x/N}.$
\end{itemize}
\et

\bp  In order to prove (a), we note that if $P\chi_{[0,N]}$ belongs to $C^\infty(\mr),$ all the derivatives at
$x=0$ vanishes. But $P$ is an entire function and therefore equal to its Taylor series, so this would imply that $P$ is identically zero, which is a contradiction.
For the proof of  (b), fix $L\in \mn$. The ``if" implication is clear, so suppose that $P  \chi_{[0, N]} \in C^{L-1}(\mr)$. We use induction to show \eqref{et-78}. First, observe that $P(0)=D
P(0)=\cdots = D^{L-1} P(0)=0.$
Since  $P(0)=\sum_{k\in \mz}c_k=0$, we have
\begin{equation*}
  P(x)=\sum_{k\neq 0 } c_k (e^{2\pi k x/N}-1).
\end{equation*}
Define $P_+$ and $P_-$ by
$$P_+(x):=\sum_{k\in \mn } c_k (e^{2\pi ik x/N}-1), \ \ \ \
P_-(x):=\sum_{k\in \mn } c_{-k} (e^{-2\pi k i x/N}-1).$$
Then we see that
\begin{eqnarray*}
   P_+(x)&=&\sum_{k\in\mn}c_k(e^{2\pi i x/N} -1)\sum_{\ell=0}^{k-1}e^{2\pi i \ell x/N} \\
         &=& e^{\pi i x/N}\sin(\pi x/N)\left(2i \sum_{k\in \mn} c_k \sum_{\ell=0}^{k-1} e^{2\pi i \ell x/N}  \right)\\
         &=:& e^{\pi i x/N}\sin(\pi x/N) \Lambda_+(x).
\end{eqnarray*}
Similarly,
\begin{eqnarray*}
P_-(x)&=&e^{\pi i x/N}\sin(\pi x/N)
\left(-2i \sum_{k\in \mn} c_{-k} \sum_{\ell=1}^{k} e^{-2\pi i \ell x/N}  \right)\\
&=:& e^{\pi i x/N}\sin(\pi x/N) \Lambda_-(x).
\end{eqnarray*}
Then we have
$$P(x)=P_+(x)+P_-(x)=e^{\pi i x/N}\sin(\pi x/N) A_1(x),$$
where $A_1(x):=\Lambda_+(x)+\Lambda_-(x)$ is an $N$-periodic function. In order to arrive at \eqref{et-78} we will now inductively  assume that,
for some $1\leq \ell \leq L-1$,
\begin{equation} \label{et-79}
P(x)= \left(e^{\pi i x/N} \sin(\pi
x/N)\right)^{\ell} A_{\ell}(x)
\end{equation}
 for an  $N$-periodic entire function $A_{\ell}$.
By the Leibnitz formula for the $\ell$th derivative of a product, we have
\begin{equation} \label{et-8}
  D^\ell P(x)= \frac{1}{(2i)^\ell}
  \sum_{k=0}^\ell {\ell \choose k}
  D^k \left(e^{2\pi i x/N}-1 \right)^\ell
  D^{\ell-k} A_{\ell} (x).
\end{equation}
Since
 $  D^k \left(e^{2\pi i x/N}-1 \right)^\ell
   = \ell(\ell-1)\cdots(\ell-k+1)\left(e^{2\pi i x/N}-1  \right)^{\ell-k}
   \left(\frac{2\pi i}{N}\right)^k,$
we have $
   D^k \left(e^{2\pi i x/N}-1 \right)^\ell \left|_{x=0} \right.
   = \ell ! \left(\frac{2\pi i}{N}\right)^\ell  \delta_{\ell, k}.$
It follows from \eqref{et-8} that \\
$D^\ell P(0) =\frac{\ell !}{(2i)^\ell}
    \left(\frac{2\pi i}{N}\right)^\ell A_{\ell} (0).$
By assumption $D^\ell P(0)=0,$ so we conclude that $A_\ell (0)=0.$
By an argument similar to the case $P(0)=0$,
we see that
$$A_{\ell} (x)
= e^{\pi i x/N} \sin (\pi x/N)  \Lambda_{\ell+1}(x)$$
for an $N$-periodic entire function $\Lambda_{\ell+1}(x)$.
This together with \eqref{et-79}
leads to
$P(x)=\left( e^{\pi i x/N} \sin (\pi x/N) \right)^{\ell+1} \Lambda_{\ell+1}(x).$
This completes the induction.
\ep

\section{Trigonometric polynomials} \label{277a}
In this section we specialize to the case of trigonometric polynomials.
Theorem \ref{et-61} characterizes the
regularity that can be obtained for $P\chi_{[0,N]}$ when $P$ is a trigonometric polynomial with a given number of terms and coefficients satisfying the condition in Lemma \ref{157d}. The
subsequent Propositions \ref{et-76} and \ref{et-77} will  show that the partition of unity condition can be combined with any finite regularity and any support size  by taking a trigonometric polynomial of sufficiently high degree.

\bt \label{et-61}

Let $K, N\in \mn$. Assume that
\begin{equation}\label{et-39}
P(x):= \sum_{k=-K}^{K} c_{k} e^{2\pi ikx/N}
\end{equation}
is a real-valued trigonometric polynomial with $c_k=\frac{1}{N}\delta_{0,k}, k\in N\mz$. Then the following hold.
\begin{itemize}
\item[\rm{(a)}] There does not exist $P$ of the form \eqref{et-39}
such that $P\chi_{[0,N]}\in C^{2K}(\mr)$;
\item[\rm{(b)}] Fix $L \in \{1,2,\cdots, 2K \}.$ Then  $P\chi_{[0,N]} \in C^{L-1}(\mr)$ if and only if \\
$P(x)= \left(e^{\pi i x/N} \sin(\pi x/N)\right)^{L} A_{L}(x)$ for a  trigonometric polynomial
\begin{equation}\label{et-69}
   A_{L}(x):= \sum_{k=-K}^{K-L} a_k e^{2\pi i k x/N}.
\end{equation}


\end{itemize}
\et

\bp (a): Assume  that there exists $P$
such that $p\in C^{2K}(\mr)$. Then $ P(0)=D P(0)=\cdots = D^{2K} P(0)=0,$ that is,
$\sum_{k=-K}^{K} k^i c_k=0,\  i=0,1,\cdots, 2K.$
This set of equations can be written in the form
$M_1\{c_k\}_{k=-K}^{K}= {\bf 0},$ where
$M_1$ is the $(2K+1)\times (2K+1)$ matrix defined by
$$M_1=\begin{pmatrix}
1 & \cdots &  1 &  1 & 1 &\cdots  &  1 \\
-K & \cdots & -1 & 0 & 1 & & N-1  \\
(-K)^2 & \cdots & (-1)^2 & 0 & 1^2 &\cdots &  K^2 \\
\vdots &  & \vdots & \vdots & \vdots &  & \vdots \\
(-K)^{2K} & \cdots & (-1)^{2K} & 0 & 1^{2K} & \cdots  & K^{2K}
\end{pmatrix}.
$$
This is a $2K+1 \times 2K+1$ Vandermonde matrix,  with rows determined by the
numbers $ z_k = -K +k, \ k=0,1,\cdots,
2K,$ and therefore invertible. Hence  the system only has
the trivial solution. This contradicts the assumption that  $c_0=1/N$.
The proof of (b) follows the lines of the proof of Theorem \ref{et-68} by keeping track of the number of terms in the trigonometric polynomials.
\ep

For any $N\in \mn$ we will now show how to construct trigonometric polynomials $P$ such that $P\chi_{[0,N]}$  satisfies the partition of unity condition and has desired regularity. We begin with the case $N=2.$

\begin{prop} \label{et-76}
Let $N=2$.   Consider a real--valued trigonometric polynomial
$Q(x)=\sum_{k} c_k e^{\pi i k x }$
with $c_{k}=\frac{1}{2}\delta_{k,0}, \ \ k\in 2\mz.$
Given  $L \in \mn,$ define a trigonometric
polynomial $P$ by
\begin{equation}\label{et-74}
   P(x):=Q^{L}(x) \sum_{k=0}^{L-1}{2L-1 \choose k} Q^{L-1 -k}(x) Q^k (x+1).
\end{equation}
Then $P  \chi_{[0,2]}$ satisfies the partition of unity property.
If $Q \, \chi_{[0,2]}\in C^{1}(\mr)$, then
$P  \chi_{[0,2]}\in C^{2 L -1}(\mr)$.

\end{prop}

\bp
Note that  $Q\chi_{[0,2]}$ satisfies the partition of unity property by Lemma \ref{157d}.
Using the Binomial Theorem, we have
\begin{eqnarray}
1=\left( Q(x) +Q(x+1)  \right)^{2L-1}
= \sum_{k=0}^{2L-1} {2L-1 \choose k} Q^{2L-1-k}(x) Q^k(x+1). \ \hspace{.3cm} \ \label{et-75}
\end{eqnarray}
Take $P$ as in \eqref{et-74}. Then
$$ P(x)= \sum_{k=0}^{2L-1}{L-1 \choose k} Q^{2L-1 -k}(x) Q^k (x+1).$$
Using the $2$-periodicity of $Q$ implies that
\begin{eqnarray*}
   P(x+1)&=& \sum_{k=0}^{L-1} {2L-1 \choose k} Q^{2L-1-k}(x+1) Q^k(x) \\
   &=& \sum_{\ell=L}^{2L-1} {2L-1 \choose \ell} Q^{\ell}(x+1) Q^{2L-1-\ell}(x).
\end{eqnarray*}
By \eqref{et-75}, we have
$P(x) +P(x+1)=1,$ so
$P\chi_{[0,2]}$ satisfies the partition of unity property, as desired. Furthermore,
if $Q \, \chi_{[0,2]} \in C^1(\mr)$, then by (b) in Theorem \ref{et-68} we know that
$Q(x)=\sin^2(\pi x/2) e^{\pi i x}A(x)$ for some 2-periodic entire function (actually a trigonometric polynomial) $A.$
Using \eqref{et-74}, it follows that
$$ P(x)=\sin^{2L}(\pi x/2) e^{\pi i xL}\tilde  A(x)$$
for a 2-periodic entire function (trigonometric polynomial) $\tilde A$.
By Theorem \ref{et-68} (b) we conclude  that $P  \chi_{[0,2]}\in C^{2L-1}(\mr)$.  \ep

In order to construct partition of unities based on functions supported on $[0,2]$ and with desired regularity, we just need to
provide a
concrete example of a trigonometric polynomials $Q$ satisfying the conditions in
Proposition \ref{et-76}: \bex \label{et-76e}
Let \bes Q(x):= \sin^2(\pi  x/2) =\left(\frac{e^{i\pi x/2}-e^{-i\pi x/2}}{2i}\right)^2 =  -\frac14 e^{\pi ix} + \frac12 - \frac14 e^{-\pi x}.\ens
Then  $Q$ has the form described in Proposition \ref{et-76},
and $Q\chi_{[0,2]}\in C^1(\mr).$
\ep \enx

Let us now consider the case $N\geq 3.$ Starting with a certain trigonometric polynomial $P_1$ such that $P_1\chi_{[0,N]}$ satisfies the partition of unity condition we provide an inductive procedure  which,  in each step, yields a new trigonometric polynomial with the partition of unity property and higher regularity.

\begin{prop} \label{et-77}
Let $N\in \mn$ with $N\geq 3$ and put $K=N/2$ if $N$ is even, $K=(N-1)/2$ if $N$ is odd.
Consider a trigonometric polynomial $A_1(x)=\sum a_k e^{2\pi i
x/N}$, and assume that for
$P_1(x):= \left(\prod_{k=0}^K \sin^2(\pi (x-k)/N) \right) A_1(x),$
$P_1\chi_{[0,N]}$
satisfies the partition of unity property. For $L\in \{2, 3, \dots\},$ let
$P_L$ be inductively defined by
\begin{eqnarray}
   && P_L(x) \label{et-71}\\
   &&:=
   \left\{
   \begin{array}{ll}
P_{L-1}(x)\left(P_{L-1}(x)  + 2 \sum_{n=1}^{K-1} P_{L-1}(x+n) + P_{L-1}(x+K)
\right),
& \text{ if $N$ is even};\\
P_{L-1}(x)\left(P_{L-1}(x)  + 2 \sum_{n=1}^{K} P_{L-1}(x+n) \right),
& \text{ if $N$ is odd}.
   \end{array}
   \right. \nonumber
\end{eqnarray}
Then the following holds:
\bei \item[(i)] $P_L$  can be factorized as
\begin{equation}\label{et-73}
P_{L}(x)= \sin^{2L}(\pi x/N) \left(\prod_{k=1}^K \sin^2(\pi (x-k)/N) \right) A_{L}(x)
\end{equation}
for some trigonometric polynomial $A_{L}(x)=\sum_k a_k^{(L)}
e^{2\pi i k x /N}$.
\item[(ii)] $P_L \, \chi_{[0,N]}$ satisfies the partition of unity
property and belongs to $C^{2L-1}(\mr)$.
\eni

\end{prop}

\bp We give the proof for the case where $N$ is even, $K=N/2;$ the case where $N$ is odd is just requires minor modifications. We will use induction. By assumption (i) and (ii) are satisfied for
$L=1.$  Now assume that for some $L\in \mn$ with $L\geq 2$, $P_{L-1}\chi_{[0,N]}$
satisfies the partition of unity property and that
\begin{equation} \label{et-70}
P_{L-1}(x)=\sin^{2(L-1)}(\pi x/N) \left(\prod_{k=1}^K \sin^2(\pi (x-k)/N)\right)A_{L-1}(x),
\end{equation}
for a trigonometric polynomial $A_{L-1}(x)=\sum_k a_k^{(L-1)}
e^{2\pi i k x /N}$. Then $P_L$  can be factorized as in \eqref{et-73} with
\begin{eqnarray*}
A_{L}(x)&:=& \sin^{2L-4}(\pi x/N) A_{L-1}^2(x) \prod_{k=1}^K \sin^2(\pi (x-k)/N)\\
&+&2 \sum_{\ell=1}^{K-1} \sin^{2L-2}(\pi (x+\ell)/N)
 A_{L-1}(x+\ell)A_{L-1}(x) \prod_{k\in\{-\ell+1,\cdots,K-\ell\}\setminus\{0\}}\sin^2(\pi (x-k)/N)\\
&+&  \sin^{2L-2}(\pi (x+K)/N)
 A_{L-1}(x+K)A_{L-1}(x) \prod_{k=1}^{K-1}\sin^2(\pi (x+k)/N).
\end{eqnarray*}
By Theorem \ref{et-68} (b) it follows that
$P_L \, \chi_{[0,N]}\in C^{2L-1}(\mr)$.
We now show $P_{L}\chi_{[0,N]}$ satisfies the partition of unity property.
From the  partition of unity property for $P_{L-1}\chi_{[0,N]}$, for $x\in [0,1]$ we have
\begin{eqnarray}
   1&=& \left( \sum_{j=0}^{N-1} P_{L-1}(x+j)  \right)^2  \nonumber\\
   &=& \sum_{j=0}^{N-1} P_{L-1}^2(x+j)+
   2 \sum_{0\leq k<n\leq N-1} P_{L-1}(x+k)P_{L-1}(x+n) \nonumber\\
   &=& \sum_{j=0}^{N-1} P_{L-1}^2(x+j)+
   2 \sum_{n=1}^{N-1} \sum_{k=0}^{N-1-n} P_{L-1}(x+k)P_{L-1}(x+k+n). \label{et-85}
\end{eqnarray}
Observe that
\begin{eqnarray}
&&\sum_{n=1}^{N-1} \sum_{k=0}^{N-1-n} P_{L-1}(x+k)P_{L-1}(x+k+n)  \nonumber\\
&&=\left(\sum_{n=1}^{K-1}\sum_{k=0}^{2K-1-n} +
\sum_{n=K+1}^{2K-1}\sum_{k=0}^{2K-1-n} \right)P_{L-1}(x+k)P_{L-1}(x+k+n)  \nonumber \\
&&+\sum_{k=0}^{K-1}P_{L-1}(x+k)P_{L-1}(x+k+K). \label{et-84}
\end{eqnarray}
Note that
\begin{eqnarray*}
  && \sum_{n=K+1}^{2K-1}\sum_{k=0}^{2K-1-n} P_{L-1}(x+k)P_{L-1}(x+k+n)  \\
  &&= \sum_{m=1}^{K-1}\sum_{k=0}^{m-1} P_{L-1}(x+k)P_{L-1}(x+k+2K-m) \ \ \ \ (m=2K-n)\\
  &&= \sum_{m=1}^{K-1}\sum_{j=2K-m}^{2K-1} P_{L-1}(x+j-2K+m )P_{L-1}(x+j) \ \ \ \ (j=k+2K-m)\\
  &&= \sum_{m=1}^{K-1}\sum_{j=2K-m}^{2K-1} P_{L-1}(x+j+m )P_{L-1}(x+j),
\end{eqnarray*}
where we used the $2K$-periodicity of  $P_{L-1}$ in the last equality.
This implies
\begin{eqnarray}
&&\left(\sum_{n=1}^{K-1}\sum_{k=0}^{2K-1-n} +
\sum_{n=K+1}^{2K-1}\sum_{k=0}^{2K-1-n} \right)P_{L-1}(x+k)P_{L-1}(x+k+n) \nonumber \\
&&=\sum_{n=1}^{K-1}\sum_{k=0}^{2K-1}  P_{L-1}(x+k)P_{L-1}(x+k+n). \label{et-82}
\end{eqnarray}
Moreover,
\begin{eqnarray*}
&&\sum_{k=0}^{K-1}P_{L-1}(x+k)P_{L-1}(x+k+K)  \\
&&=\sum_{j=K}^{2K-1}P_{L-1}(x+j-K)P_{L-1}(x+j)  \ \ \ \ (j=k+K) \\
&&=\sum_{j=K}^{2K-1}P_{L-1}(x+j+K)P_{L-1}(x+j).
\end{eqnarray*}
Then we obtain
\begin{eqnarray*}
\sum_{k=0}^{K-1}P_{L-1}(x+k)P_{L-1}(x+k+K)
=\frac{1}{2}\sum_{k=0}^{2K-1}P_{L-1}(x+k)P_{L-1}(x+k+K).
\end{eqnarray*}
Putting this and \eqref{et-82} into the right-hand side of \eqref{et-84}, we have
\begin{eqnarray*}
&&\sum_{n=1}^{2K-1} \sum_{k=0}^{2K-1-n} P_{L-1}(x+k)P_{L-1}(x+k+n)\\
&&= \sum_{k=0}^{2K-1} P_{L-1}(x+k)
\left( \sum_{n=1}^{K-1} P_{L-1}(x+k+n) +\frac{1}{2}P_{L-1}(x+k+K) \right).
\end{eqnarray*}
Combining this with \eqref{et-85} yields
\begin{eqnarray*}
   1&=& \sum_{j=0}^{2K-1}P_{L-1}(x+j)
   \left(P_{L-1}(x+j)+ 2\sum_{n=1}^{K-1} P_{L-1}(x+j+n) +P_{L-1}(x+j+K)  \right)\\
   &=& \sum_{j=0}^{2K-1} P_L(x+j).
\end{eqnarray*}
Hence  $P_{L}\chi_{[0,N]}$ satisfies the partition of unity property, as desired.   \ep

Again, in order to construct partition of unities based on functions supported on $[0,N]$ and with desired regularity, we just need to
provide a suitable polynomial $P_1$ in Proposition \ref{et-77}:

\bex Given $N\in \{3,4, \dots,\},$ put $K=N/2$ if $N$ is even, $K=(N-1)/2$ if $N$ is odd. Let
\begin{equation*}
P_1(x):= \left(\sum_{n=K+1}^{N-1}
\prod_{k=0}^{K} \sin^2(\pi (n-k)/N) \right)^{-1}
\prod_{k=0}^K \sin^2(\pi (x-k)/N).
\end{equation*}
A direct calculation shows that $P_1$ has the form
$$P_1(x)=\sum_{k=-K-1}^{K+1} c_k e^{2\pi i k x/N}.$$
Note that since  $N\geq 3$, we have $K+1<N$.
Thus
\begin{equation*}
1+ e^{2\pi i k/N}
+ \cdots + \left( e^{2\pi i k/N}\right)^{N-1}=\left\{
\begin{array}{ll}
   N, & k=0 \\
   0, & k\in\{\pm 1, \cdots, \pm (K+1) \},
\end{array}
\right.
\end{equation*}
and
\begin{eqnarray*}
&& P_1(x)+P_1(x+1)+\cdots+P_1(x+N-1) \\
&&=\sum_{k=-K-1}^{K+1} c_k\,\left[1+ e^{2\pi i k/N}
+ \cdots + \left( e^{2\pi i k/N}\right)^{N-1} \right] \, e^{2\pi ikx/N}=c_0 N.
\end{eqnarray*} We must now determine the constant $c_0N.$
Taking $x=0$ we see that
\bes & \ & P_1(0)+P_1(1)+P_1(2)+\cdots+P_1(N-1)  \\ &=&\left(\sum_{n=K+1}^{N-1}
\prod_{k=0}^{K} \sin^2(\pi (n-k)/N) \right)^{-1}
\sum_{n=0}^{N-1}\prod_{k=0}^K \sin^2(\pi (n-k)/N)\\ &=&\left(\sum_{n=K+1}^{N-1}
\prod_{k=0}^{K} \sin^2(\pi (n-k)/N) \right)^{-1}
\sum_{n=K+1}^{N-1}\prod_{k=0}^K \sin^2(\pi (n-k)/N)
=1.\ens
It follows that $P_1 \chi_{[0,N]}$ satisfies the partition of unity property.
\ep
\enx


\section{Applications to Gabor frames} \label{157g}

To recapitulate, the results in Section \ref{277a} shows that for each $N,L\in \mn$ we can construct a trigonometric polynomial $P$ such that the restriction to $\mr$ is real valued, $P\chi_{[0,N]}\in C^L(\mr),$ and $P\chi_{[0,N]}$ satisfies the partition of unity condition.  The purpose of this section is to provide new constructions of pairs of dual Gabor frames based on these functions.  We will
restrict our presentation to the case where the translation parameter is $a=1,$ but  via a scaling, other choices of $a$ are possible if we change the length of the support accordingly.
 As starting point, let us consider the following result from \cite{CR}.

\bpr \label{128a} Let $N\in \mathbb{N}$. Let $g\in L^2(\mathbb{R})$ be a real-valued
bounded function for which $supp \ g \subseteq [0,N]$ and
$\sum_{n\in \mz} g(x-n)=1. $ Let $b\in
]0,\frac1{2N-1}]$. Define $h\in L^2(\mathbb{R})$ by \bes
h(x)=\sum_{n=-N+1}^{N-1} a_n g(x+n), \ens where
$a_0=b, \  \ a_n + a_{-n}=2b, \ n=1,2,\cdots, N-1.$ Then $g$
and $h$ generate dual frames $\{E_{mb} T_n g\}_{m,n\in\mathbb{Z}}$
and $\{E_{mb} T_n h\}_{m,n\in\mathbb{Z}}$ for $L^2(\mathbb{R}).$\epr

The constructions in Proposition \ref{et-76} and Proposition \ref{et-77} have exactly the properties required in Proposition \ref{128a}, combined with desired regularity. We will formulate the corresponding construction of a pair of dual pair of frames based on the case $N=2,$ and leave the formulation for the case $N\ge 3$
to the reader.

\bc \label{128b} Take the trigonometric polynomial  $Q$ as in
Proposition \ref{et-76} and let, for $L\in \mn,$
\bes\label{et-74a}
   P(x):=Q^{L}(x) \sum_{k=0}^{L-1}{2L-1 \choose k} Q^{L-1 -k}(x) Q^k (x+1).
\ens Let $b\in ]0, 1/3]$ and choose coefficients $a_n$ such that
$a_0=b, \  \ a_1 + a_{-1}=2b.$ Then the functions
\begin{equation} \label{et-86}
g(x):=(P\chi_{[0,2]})(x) \ \mbox{ and} \ \ h(x):=\sum_{n=-1}^{1}
a_n (P\chi_{[0,2]})(x+n)
\end{equation}
 belong to $C^{2L-1}(\mr)$ and
generate dual Gabor frames $\{E_{mb} T_n g\}_{m,n\in\mathbb{Z}}$
and $\{E_{mb} T_n h\}_{m,n\in\mathbb{Z}}$ for
$L^2(\mathbb{R}).$\ec

\begin{ex} \label{et-87}
Let $L=2, b\in ]0,1/3],$ and $Q(x):= \sin^2(\pi  x/2).$
Define
$$   P(x)=\sin^4(\pi x/2) \sum_{k=0}^{1}{3 \choose k}
\sin^{2(1 -k)}(\pi x/2) \sin^{2k}(\pi (x+1)/2).
$$
Then
$g$ and $h$
defined as in \eqref{et-86}
 belong to $C^{3}(\mr)$ and
generate dual Gabor frames $\{E_{mb} T_n g\}_{m,n\in\mathbb{Z}}$
and $\{E_{mb} T_n h\}_{m,n\in\mathbb{Z}}$ for
$L^2(\mathbb{R}).$
On Figure \ref{et-f1}, we plot $g$ and $h$ for the choice  $a_{-1}=a_0=a_1=1/3.$  \ep

\begin{figure}
\centerline{
\subfigure[]{\includegraphics[width=2.7in,height=2in]{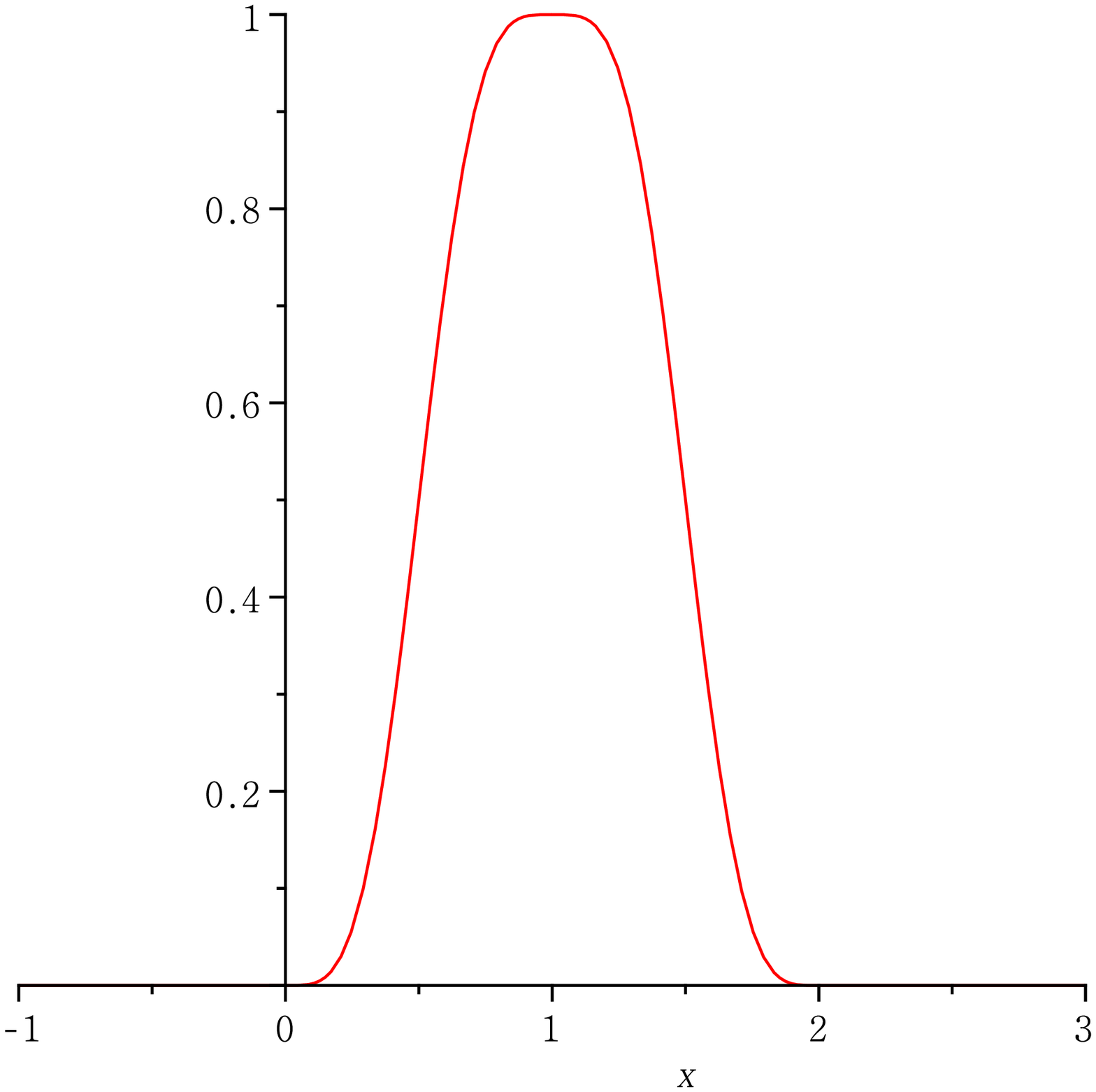}}\hfil
\subfigure[]{\includegraphics[width=2.7in,height=2in]{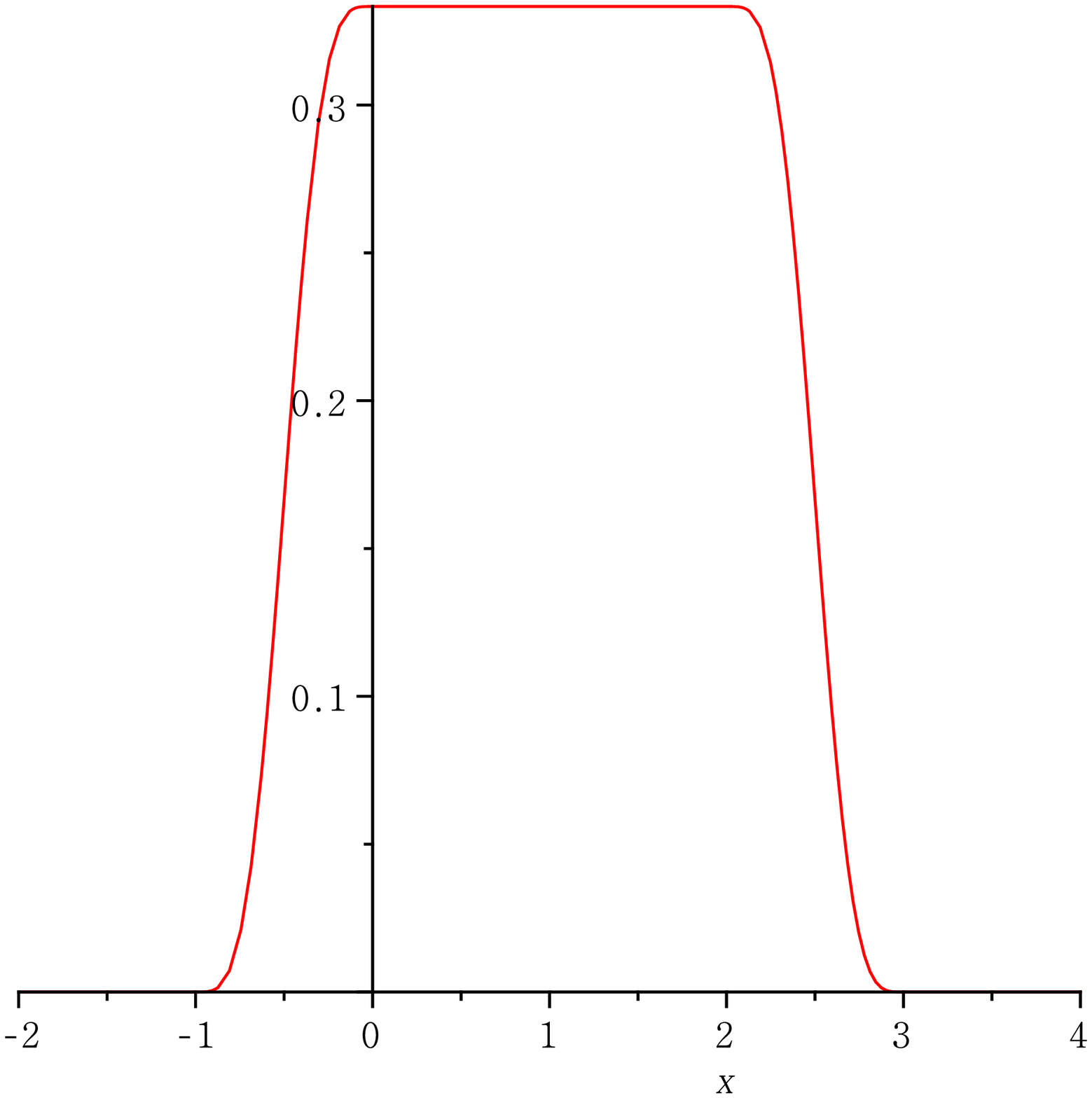}}}\hfil
\caption{ Plots of the dual windows in Example \ref{et-87}:   (a)
$g$;  (b) $h$ with $a_{-1}=a_0=a_1=1/3.$
\label{et-f1}}
\end{figure}

\end{ex}

Compared with the other results in the literature, Corollary \ref{128b} has the advantage that desired regularity of the frame generators does not make the support size grow and the redundancy increase. In the classical application of Proposition \ref{128a} where the function $g$ is a B-spline $B_N$ for some $N\in \mn,$ high regularity can only be obtained by considering large values for $N,$ which leads to
functions with large support and corresponding small values for the parameter $b.$ Since the redundancy of a Gabor frame $\mts$ is measured by the number $1/(ab)=1/b,$ these constructions have high redundancy.
On the other hand the frames in Corollary \ref{128b} are generated by windows that are supported on $[0,2],$ the dual windows are supported on $[-1, 3]$ regardless of the desired regularity, and by taking $b=1/3$ the redundancy is just 3.

In the setting discussed here, we can even avoid to choose a dual window with larger support than the given window. Hereby we can enlarge the range of the parameter $b,$ and provide constructions with redundancy 2. Note that also the construction by Laugesen in \cite{L} keeps the support size.

\bc \label{et-80}
Let $L_1, L_2\in \mn$, and fix $b\in ]0, \frac{1}{2}]$.
Take $P(x):=\sin^2(\pi x/2)$.
Define $$g(x)=\sin^{2L_1}(\pi x /2)\chi_{[0,2]}(x)$$
and
$$h(x)=b \sin^{2L_2}(\pi x /2)
\left( \sum_{k=0}^{L_1+L_2-1}{2L_1+2L_2-1 \choose k} P^{L_1+L_2-1
-k}(x) P^k (x+1) \right)\chi_{[0,2]}(x).$$ Then $g \in
C^{2L_1-1}(\mr)$, $h \in C^{2L_2-1}(\mr)$, and the functions
 $\{E_{mb}T_ng\}_{m,n\in \mz}$ and $\{E_{mb}T_n h\}_{m,n\in \mz}$ form
a pair of dual frames.
\ec
\bp Using Proposition \ref{et-76} together with Example \ref{et-76e}, it follows that the function $gh$ satisfies the condition \eqref{gfs} for $n=0.$ The choice of $b$ and the support sizes for $g$ and $h$ shows that \eqref{gfs} holds for $n\neq 0$ as well.\ep

\begin{figure}
\centerline{
\subfigure[]{\includegraphics[width=2.7in,height=2.2in]{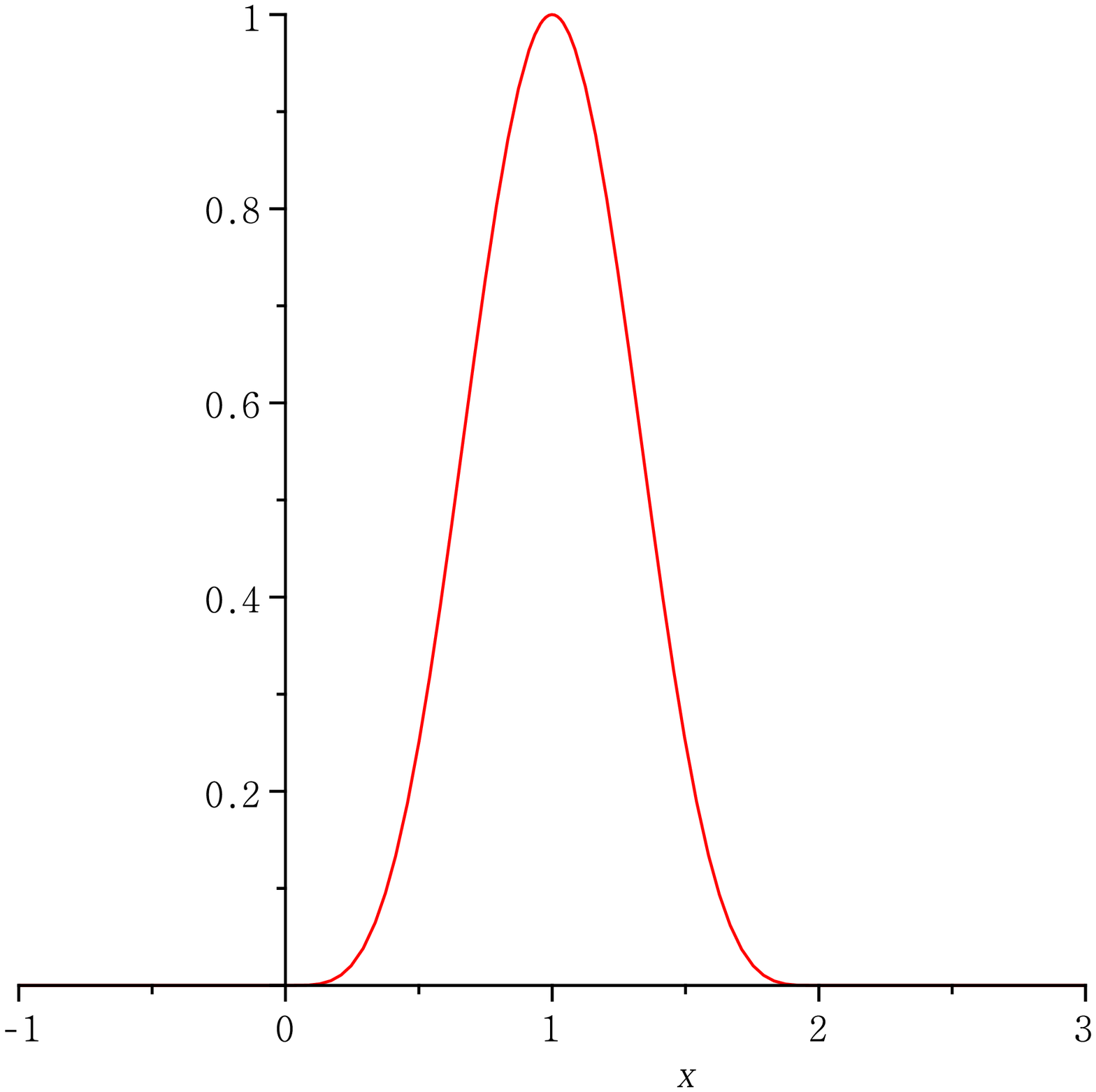}}\hfil
\subfigure[]{\includegraphics[width=2.7in,height=2.2in]{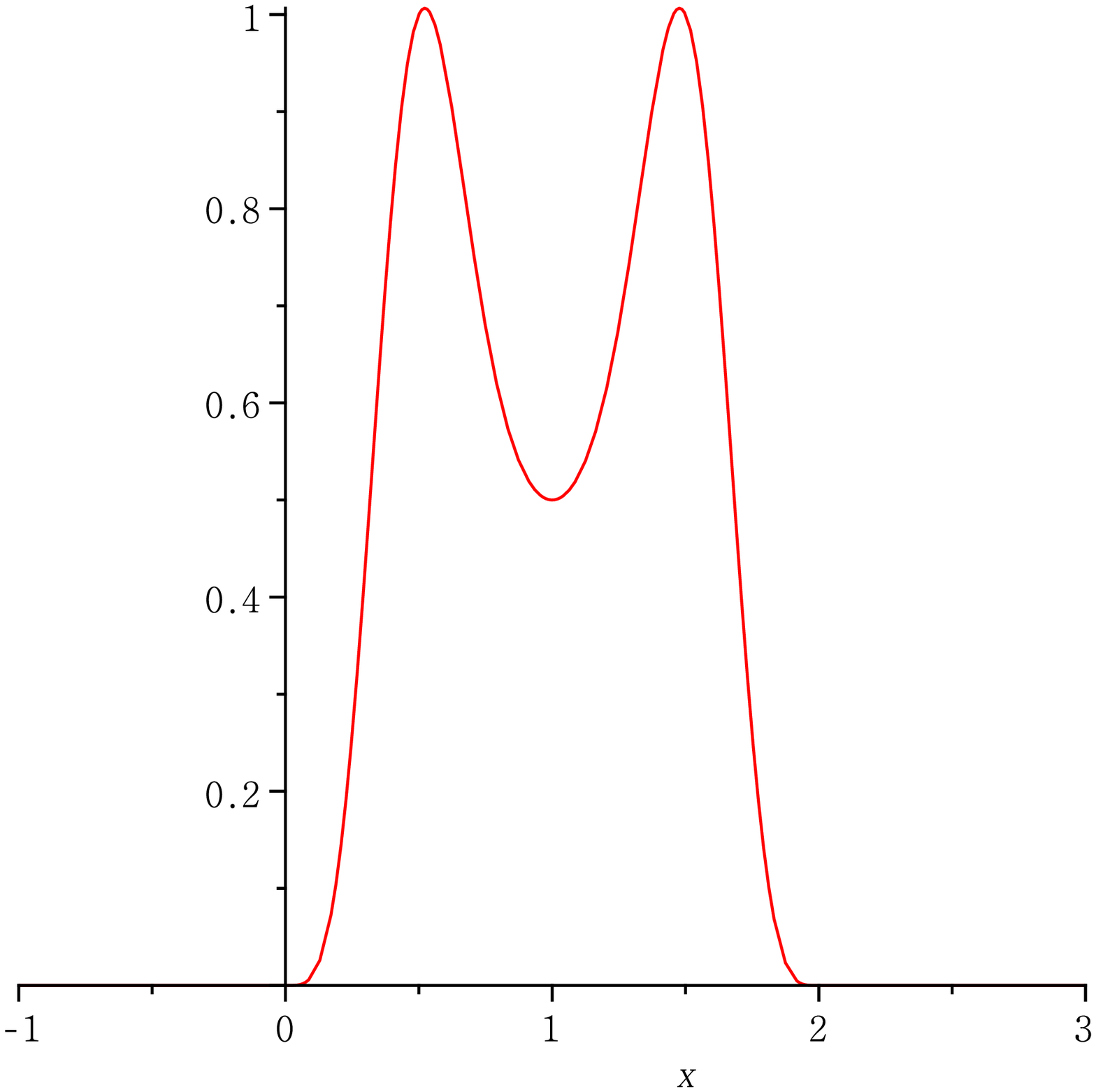}}}\hfil
\caption{ The windows $g$ and $h$ in Corollary \ref{et-80} for $L_1=L_2=2, b=1/2$
\label{et-f2}}
\end{figure}

\noindent Figure \ref{et-f2} shows the windows $g$ and $h$ for $L_1=L_2=2, b=1/2.$

\newpage
Up to a certain regularity an even simpler construction is possible:

\bc \label{237b} Let $N\in \mn$ and fix $b\in
]0, \frac{1}{N}].$ Then the following hold:
  \bei \item[(i)] For   $L_1, L_2 \in \mn$ with $1\leq L_1 +L_2 \leq N-1$, the functions
$$g(x)=\sin^{2 L_1}(\pi x /N) \chi_{[0,N]}(x), \hspace{.3cm} \
h(x)=\frac{b 4^{L_1+L_2}}{N {2(L_1+L_2) \choose L_1+L_2}} \sin^{2 L_2}(\pi x /N)\chi_{[0,N]}(x)$$
belong to $ C^{2L_1-1}(\mr) $ and $C^{2L_2-1}(\mr) $ respectively, and  generate a pair of dual frames
$\{E_{mb}T_ng\}_{m,n\in \mz}$ and $\{E_{mb}T_nh\}_{m,n\in \mz}$.
\item[(ii)] For any positive integer $L \le N-1$ the function
$g(x):= \sqrt{ \frac{b4^L}{N {2L \choose L}}} \sin^L(\pi x/N)\chi_{[0,N]}(x)$
belongs to $C^{L-1}(\mr)$ and generates a tight Gabor frame $\{E_{mb}T_ng\}_{m,n\in \mz}$
\eni
\ec
\bp
Let
$$P(x):=\frac{b 4^{L_1+L_2}}{N {2(L_1+L_2) \choose L_1+L_2}} \sin^{2(L_1+L_2)}(\pi x /N).$$
By Theorem  \ref{et-61}, to prove (i) it suffices to show that the Fourier coefficients $c_k$
in the expansion $P(x)=\sum_{k\in \mz}c_k e^{2\pi i k x/N}$
satisfy $c_k=\frac{1}{N}\delta_{0,k},\ k\in N\mz$.
Via the Binomial formula,
\begin{eqnarray*}
  \sin^{2(L_1+L_2)}(\pi x/N)&=&
  \left( \frac{e^{\pi  i x /N} -e^{-\pi i x/N}}{2i} \right)^{2(L_1+L_2)}\\
 &=& \frac{1}{4^{L_1+L_2}} \sum_{k=-L_1-L_2}^{L_1+L_2} (-1)^k {2(L_1+L_2)\choose L_1+L_2+k}
 e^{2\pi i k x/N}.
 \end{eqnarray*}
The result in (ii) follows by taking $L=L_1+L_2.$\ep

Note that construction of Gabor frames based on trigonometric
polynomials appear at other places in the literature. In \cite{DGM},
Daubechies, Grossmann and Meyer construct a tight Gabor frame
based on the function $g(x)= \cos(x) \ \chi_{[-\pi/2, \pi/2]}(x),$ which is just a shifted and scaled version of the function $\sin(\pi x/2)\chi_{[0,2]}.$
Also, in \cite{CM}, the authors consider frames generated by
functions of the form $g_k(x)= \sin^k(\pi x/3) \, \chi_{[0,3]}(x)$
for parameters $k \in \mn.$ Interestingly, the results in
\cite{CM} show that $g_k$ generates a frame for all $b\in ]0,
1/3]$ and all $k \in \mn;$ but only for $k <6$ there is a dual
Gabor frame $\{E_{mb}T_nh\}_{m,n\in \mz}$  for a function of the
form $h(x)= a_0 g_k(x)+ a_1 g_k(x+1)+ a_2g_k(x+2).$ Since
$g_k$ corresponds with our setup for $N=3,$ this result is in
accordance with the limitation on the possible parameters $L_1, L_2$ in
Corollary \ref{237b}.

Note that the Corollaries  \ref{et-80} and  \ref{237b} are restricted to the case $b\le 1/N.$
For constructions of dual frame pairs generated by entire functions, we can actually show that this condition is necessary:


\bpr Let $N\in \mn$ and let $G$ and $H$
be real-valued, $N$-periodic entire functions.
If $g:=G \chi_{[0,N]}$ and $h:=H \chi_{[0,N]}$
generate  dual  frames
$\{E_{mb}T_ng\}_{m,n\in \mz}$ and $\{E_{mb}T_nh\}_{m,n\in \mz},$ then
$0< b\leq 1/N.$
\epr
\bp
To get a contradiction, assume that $\{E_{mb}T_ng\}_{m,n\in \mz}$ and $\{E_{mb}T_nh\}_{m,n\in \mz},$ are dual frames and that $1/N< b <1$.
There exists a unique $n\in\{1,2,\cdots, N-1\}$ such that
$n\leq 1/b <n+1.$
From the duality conditions, we have
\begin{equation} \label{et-62}
\sum_{j\in \mz}
G(x+j) \chi_{[0,N]}(x+j) H(x+j+1/b)\chi_{[0,N]}(x+j+1/b)= 0.
\end{equation}
For $0 < x < n+1 -1/b(<1)$,  we have
$n\leq \frac{1}{b}<x+\frac{1}{b}<n+1.$
From \eqref{et-62}, for $x\in ]0,n+1-1/b[,$
\begin{equation}\label{et-63}
\sum_{j=0}^{N-n-1}G(x+j)H(x+j+1/b)=0.
\end{equation}
Since the finite sum \eqref{et-63} is an entire function,
\eqref{et-63} holds for all $x\in \mr$.
Due to the assumption that $G$ and $H$ are $N$-periodic and entire, we can write
them as absolutely convergent Fourier series,
$G(x)=\sum_{k\in\mz}g_k e^{2\pi i k x/N}$ and $H(x)=\sum_{k\in\mz}h_k e^{2\pi i k x/N}.$ The absolute convergence assures that the change of summation below is legitimate in the following expansion
of \eqref{et-63}:
For $ x\in \mr$,
\begin{eqnarray*}
&&0=  \sum_{j=0}^{N-n-1} \sum_{k\in \mz} g_k e^{2\pi i k (x+j)/N}
   \sum_{\ell\in\mz} h_\ell e^{2\pi i \ell (x+j +1/b)/N} \\
&&= \sum_{m\in\mz}\left(\sum_{k+\ell=m}
g_k h_\ell e^{2\pi i\ell/(bN)} \right)
\left(\sum_{j=0}^{N-n-1}e^{2\pi i j m/N}  \right) e^{2\pi i x m /N}.
\end{eqnarray*}
Note that
\begin{eqnarray*}
\sum_{j=0}^{N-n-1} e^{2\pi i j m/N}
= \sum_{j=0}^{N-n-1} \left(e^{2\pi i m /N} \right)^j
=\left\{
\begin{array}{ll}
N-n, & m\in N\mz \\
\frac{1-e^{2\pi i m(N-n)/N}}{1-e^{2\pi i m/N}}, & m\neq N\mz
\end{array}
\right.
\neq 0.
\end{eqnarray*}
Hence, we have
$$\sum_{k+\ell=m} g_k h_\ell e^{2\pi i\ell/(bN)}=0, \ \forall m\in \mz.$$
Since
\begin{eqnarray*}
&&\left(\sum_{k\in\mz} g_k e^{2\pi i k x/N}  \right)
\left(\sum_{\ell\in\mz} h_\ell e^{2\pi i \ell/(bN)} e^{2\pi i \ell x/N}  \right)\\
&&=\sum_{m\in \mz} \left( \sum_{k+\ell=m} g_k h_\ell e^{2\pi i\ell/(bN)} \right)e^{2\pi i m x/N}=0,
\end{eqnarray*}
we have that either all $g_k=0$ or all $h_k=0$, $i.e.$, either
$G\equiv 0$ or $H\equiv 0$, which is a contradiction. \ep

{\bf \vspace{.1in}

\noindent Ole Christensen\\
Department of Applied Mathematics and Computer Science\\
Technical University of Denmark,
Building 303,
2800 Lyngby  \\
Denmark \\
 Email: ochr@dtu.dk

\vspace{.1in}\noindent Hong Oh Kim \\
Department of Mathematical Sciences, KAIST\\
373-1,
Guseong-dong, Yuseong-gu, Daejeon, 305-701\\
Republic of Korea\\
Email: kimhong@kaist.edu

\vspace{.1in} \noindent Rae Young Kim \\
Department of Mathematics,
Yeungnam University\\
214-1, Dae-dong, Gyeongsan-si, Gyeongsangbuk-do, 712-749\\
Republic of Korea\\
Email:  rykim@ynu.ac.kr
}

\end{document}